\documentstyle[amssymb]{amsart}

\newcommand{\ann}{\mathrm{ann}}

\newcommand{\gldim}{\mathrm{gl.dim}}

\newcommand{\Ker}{\mathrm{Ker}\;}

\renewcommand{\th}{\theta}
\newcommand{\al}{\alpha}

\newcommand{\be}{\beta}
\newcommand{\gam}{\gamma}

\newcommand{\iso}{\cong}

\newtheorem*{thm*}{Theorem}

\newtheorem*{lem*}{Lemma}

\newtheorem*{prop*}{Proposition}

\newtheorem*{cor*}{Corollary}

\newtheorem*{defn*}{Definition}

\newtheorem*{notadefn*}{Notation and Definition}

\newtheorem*{nota*}{Notation}

\newtheorem*{note*}{Remark}

\newtheorem*{ex*}{Example}

\newenvironment{proof}{\par\noindent{\bf Proof.}}{$\square$\par\bigskip}

\begin{document}

\title[Down-up Algebras and their representation theory]
      {\Large \bf Down-up Algebras and their representation theory}
\author[Paula A.A.B. Carvalho and Ian M. Musson]
       {Paula A.A.B. Carvalho $^1$ \\ Ian M. Musson $^2$}
\address{Departamento de Matem\'atica Pura\\
         Faculdade de Ci\^encias\\
         Universidade do Porto\\
         P. Gomes Teixeira\\
         4050 Porto\\
         PORTUGAL}
\email{pbcarval@@fc.up.pt}
\address{University of Wisconsin-Milwaukee\\
         PO Box 413\\
         Milwaukee, WI 53201 U.S.A.}
\email{musson@@csd.uwm.edu}
\thanks{$^1$ The author would like to express her gratitude to the department of Mathematical Sciences
of the University of Wisconsin-Milwaukee and especially to
Professors A.D. Bell and I. M. Musson. The research was supported
by PRAXIS XXI, grant BPD/11821/97 and partially by Praxis XXI,
Project "\'Algebra e Matem\'aticas Discretas". The last version
was done within the activities of Centro de Matem\'atica da
Universidade do Porto.}
\thanks{$^2$ Research partially supported by NSF grant DMS-9801151.}
\date\today
\maketitle

\begin{abstract}
A class of algebras called down-up algebras was introduced by G. Benkart and 
T. Roby
\cite{BenkartRoby}. We classify the finite dimensional simple modules over 
Noetherian down-up algebras
and show that in some cases every finite dimensional module is semisimple.
We also study the question of when two down-up algebras are isomorphic.
\end{abstract}

\section{Introduction}

Given a field $K$ and $\al,\be,\gam$ arbitrary elements of $K$, the associative algebra
$A=A(\al,\be,\gam)$ over $K$ with generators $d,u$ and defining relations
               $$(R1)\qquad d^2u=\al dud+\be ud^2+\gam d$$
               $$(R2)\qquad du^2=\al udu+\be u^2d+\gam u$$
is a down-up algebra.

In \cite{KMP} it is shown that $A(\al,\be,\gam)$ is Noetherian if and only if $\be\neq 0$.
Down-up algebras are also studied in \cite{Benkart}, \cite{BenkartRoby}, \cite{K} and \cite{Zhao}.
In this paper we study the representation theory and the isomorphism problem for Noetherian down-up
algebras.

By \cite[2.2]{KMP} if $\be\neq 0$ then $A=A(\al,\be,\gam)$ embeds in a skew group ring
$S=R[z,z^{-1};\sigma]$. Here $R=K[x,y]$, $\sigma$ is the automorphism of $R$ defined by
$\sigma(x)=y$ and $\sigma(y)=\al y+\be x+\gam$.
As a right $R$-module, $S$ is free on the basis $\{z^n|n\in{\Bbb Z}\}$
and the multiplication in $S$ is defined by $rz=z\sigma(r)$.
The embedding $\th:A\longrightarrow S$ is given by $\th(d)=z^{-1}$,
$\th (u)=xz$, so that $\th (ud)=x$ and $\th (du)=y$.

Techniques involving skew group rings play an important role in this paper. We show
that $S$ as defined above is isomorphic to $A_d$ the localization of $A$ at
$\{d^n|n\geq 0\}$. Similarly the localization $A_u$ is a skew group ring.

In section 2 we determine the finite dimensional modules over down-up algebras.
Here we use techniques of D. Jordan \cite{Jordan93}.

We say that a left $A$-module $M$ is $d$-torsion (resp. $u$-torsion) if
$A_d\otimes_AM=0$ (resp. $A_u\otimes_AM=0$). If $A$ is a down-up algebra arising from a
finite poset then the defining representation of $A$ is both $d$ and $u$ torsion
(see \cite{BenkartRoby} for background and some examples).
Clearly skew group ring methods can tell us nothing about such modules.
In section 3 we study finitely generated modules which are both
$d$-torsion and $u$-torsion.
In particular we obtain necessary and sufficient conditions for all
such modules to be finite dimensional.

The question of when two down-up algebras are isomorphic was
raised by G. Benkart and T. Roby in \cite{BenkartRoby}. They
divided down-up into four types such that no two algebras of
different types can be isomorphic. In section 4 we solve the
isomorphism problem for Noetherian down-up algebras in three of
their cases and in the last case for fields of characteristic $0$.

A particularly interesting class of down-up algebras arises in the
following way. For $\eta\neq 0$, let $A_{\eta}$ be the algebra
with generators $h,e,f$ and relations $$he-eh=e,$$ $$hf-fh=-f,$$
$$ef-\eta fe=h.$$ Then $A_{\eta}$ is isomorphic to a down-up
algebra $A(1+\eta,-\eta,1)$ and conversely any down-up algebra
$A(\al,\be,\gam)$ with $\be\neq 0\neq \gam$ and $\al+\be=1$ has
the above form. For this and other reasons we write $\eta=-\be$
throughout this paper. Note that $A_1\iso U(sl_2)$ and $A_{-1}\iso
U(osp(1,2))$, the enveloping algebra of $sl_2$ and $osp(1,2)$,
respectively. In section 5 we study the representation theory of
the algebras $A_{\eta}$, for $\eta\neq 0$ in detail.  In
particular we give a necessary and sufficient conditions for every
finite dimensional $A_{\eta}$-module to be semisimple.

Let ${\Bbb N}$ denote the set of positive integers and ${\Bbb N}_0={\Bbb N}\cup\{0\}$.

Throughout this paper we will assume that $K$ is an algebraically closed field.

\subsection{\label{grading}}
Note that if we define $deg(d)=1$ and $deg(u)=-1$, then relations
$(R1)$ and $(R2)$ are homogeneous of degree $1,-1$. It follows
then that $A$ is a ${\Bbb Z}$-graded ring $$A=\oplus_{n\in{\Bbb
Z}} A(n).$$ Moreover if $\be\neq 0$ then using the embedding of
$A$ into $S$ easily one sees that $A(0)=R=K[ud, du]$,
$A(n)=d^nA(0)=A(0)d^n$ if $n\geq 0$ and
$A(n)=u^{-n}A(0)=A(0)u^{-n}$ if $n\leq 0$.

\subsection{\label{I1}}Unless otherwise stated, we will assume $\be\neq 0$.
Also write $\eta=-\be$.

\begin{lem*}
Let $A=A(\al,\be,\gam)$ be a down-up algebra. The sets $\{d^n|n\in
{\Bbb N}_0\}$ and $\{u^n|n\in {\Bbb N}_0\}$ are Ore sets in $A$.
\end{lem*}
\begin{proof}
Clearly $d=z^{-1}$ is a unit in $S$. Let $$B=R[z^{-1};
\sigma]=\{\sum_{i=0}^m z^{-i}r_i|r_i\in R\}$$ so that $B\subseteq
A \subseteq S=R[z,z^{-1};\sigma]$. Also any element of $S$ has the
form $d^{-n}a=\sigma^{-n}(a)d^{-n}$, where $n\geq 0$ and $a\in
B\subseteq A$. Thus $\{d^n|n\in{\Bbb N}_0\}$ is (left and right)
Ore in $A$ and $S=A_d$. By \cite[\S2]{BenkartRoby} there is an
antiautomorphism of $A$ interchanging $u$ and $d$ and it follows
that $\{u^n|n\in{\Bbb N}_0\}$ is Ore in $A$.
\end{proof}

\subsection{\label{I2}} {\bf Lemma.\;}\;{\em Any unit in $A$ belongs to $K^*$.}
\begin{proof}
It is well known that any unit in $S$ has the form $az^n$ with $a\in K^*$ and $n\in {\Bbb Z}$,
\cite[Proposition VI.1.6]{Sehgal}.
The result follows since $z^n$ is not a unit in $A$ unless $n=0$.
\end{proof}

\subsection{\label{S1}}
Let $R=K[x,y]$ and $\sigma$ defined as before. Let $f(\lambda)=\lambda^2-\al\lambda-\be$
and $\lambda_1, \lambda_2$ its roots.

Note that $\sigma$ stabilizes the subspace $W$ of $R$ spanned by
$1$, $x$ and $y$. It is useful to find $w_1$, $w_2$ in $R$ such
that $1$, $w_1$, $w_2$ is a basis for $W$ and the matrix of
$\sigma$ with respect to this basis is in Jordan canonical form.
We can take $w_1$, $w_2$ as follows:

Case 1: $\al^2+4\be\neq0$ and $\al+\be\neq1$. Then $\lambda_1$, $\lambda_2$ are distinct
and both different from 1. We set
$$w_i=\be(\lambda_i-1)x+\lambda_i(\lambda_i-1)y+\gam\lambda_i,$$
for all $i\in\{1,2\}$. Then $\sigma(w_i)=\lambda_iw_i$ for all $i\in\{1,2\}$.

Case 2: $\al^2+4\be\neq0$ and $\al+\be=1$. In this case $\al\ne
0,2$ and $f(\lambda)$ has roots $1$ and $\eta$. Set $$w_1=\be
x+y$$ $$w_2=-x+y+\gam(\al-2)^{-1}.$$ Then $\sigma(w_1)=w_1+\gam$
and $\sigma(w_2)=\eta w_2$.

Case 3: $\al^2+4\be=0$ and $\al+\be\neq1$. Then $f(\lambda)$ has a multiple root
$\al/2$. Set
$$w_1=(2\be+\al)x+(\al-2)y+2\gam$$
$$w_2=2y-2x.$$
Then $\sigma(w_1)=(\al/2)w_1$ and $\sigma(w_2)=(\al/2)w_2+w_1$.

Case 4: $\al^2+4\be=0$ and $\al+\be=1$. Then $(\al,\be)=(2,-1)$
and $1$ is a multiple root of $f(\lambda)$. Set
$$w_1=-x+y+\gam$$
$$w_2=y.$$ Then
$\sigma(w_1)=w_1+\gam$ and
$\sigma(w_2)=w_2+w_1$.

\section{Finite Dimensional Simple Modules over Down-up Algebras}

\subsection{\label{S2}} We need a construction which is similar to the one
given by D. Jordan in \cite[3.1]{Jordan93}.
Suppose $P$ is a maximal ideal of $R$ such that $\sigma^n(P)=P$ for some $n\in{\Bbb N}$ and
suppose $n$ is minimal with this property. Set
$$M_P=\bigoplus_{i=0}^{n-1} R/(\sigma^{i}(P)).$$
We can make $M_P$ a left $S$-module by defining for each $i\in\{0,\ldots,n\}$
 and $r\in R$,
$$z.(r+\sigma^{i}(P))=\sigma^{-1}(r)+\sigma^{i-1}(P).$$

\begin{lem*}
Every finite dimensional $d$-torsion free simple left $A$-module is
isomorphic to $M_P$
for some maximal ideal $P$ of $R$.
\end{lem*}
\begin{proof}
Since $A_d\iso S$, $M_P$ is a torsion free left $A$-module.
It is easy to show that $M_P$ is simple.

Conversely, assume that $M$ is a finite dimensional $d$-torsion
free simple left $A$-module. From the Euclidean algorithm and the
fact that $M$ is $d$-torsion free it is easy to conclude that
$dM=M$. Hence  we can identify $M$ with the $A_d$-module
$A_d\otimes_A M$. As M is finite dimensional, M has a finite
composition series as a $R=K[x,y]$-module with composition factors
isomorphic to $R/P_i$, for some finite number of distinct maximal
ideals $P_1,\ldots,P_n$ of $R$. As an $R$-module, we have that
$$M=\oplus_{i=1}^n M(P_i)$$ where $M(P_i)=\{m\in M|
P^{k_i}_im=0\}$ for some $k_i\in{\Bbb N}$ is an $R$-submodule of
$M$. As for each $i\in\{1,\ldots,m\}$,
$zM(P_i)=M(\sigma^{-1}(P_i))$, we conclude that the maximal ideals
$P_1,\ldots,P_n$ are all in a single orbit, and that this orbit is
finite.

If $M_i=\{m\in M|P_im=0\}$, then $\oplus_{i=0}^{n}M_i$ is an
$A_d$-submodule of $M$. Hence we have $M=\oplus_{i=0}^{n}M_i$, and
the result follows.
\end{proof}

\subsection{\label{S3'}} 
We investigate when the $\sigma$-orbit of a given maximal ideal is
finite. The proof of the following lemma will be omitted as it is
straightforward.

\begin{lem*}
Let $A=A(\al,\be,\gam)$ be a down-up algebra and $P=(w_1-a_1,w_2-a_2)$ for
$(a_1,a_2)\in K^2$. Then there is $n>0$ such that $\sigma^n(P)=P$ if and only if
one of the following holds:
\begin{itemize}
 \item[i)] In  case 1, $(\lambda_i^n-1)a_i=0$, for $i=1,2$;
 \item[ii)] In case 2, $n\gam=(\eta^n-1)a_2=0$;
 \item[iii)] In case 3, $a_1=[(\al/2)^n-1]a_2=0$;
 \item[iv)] In case 4,  $n\gam=a_1=0.$
\end{itemize}
\end{lem*}

Following the notation in \cite{BenkartRoby}, we denote by $V(\lambda)$ the Verma module
of highest weight $\lambda$. Let $\lambda_{-1}=0, \lambda_0=\lambda$ and define for each
$n\in {\Bbb N}$, $\lambda_n=\al\lambda_{n-1}+\be\lambda_{n-2}+\gam$. The Verma module
$V(\lambda)$ has basis $\{v_n|n\in {\Bbb N}_0\}$. The action of $A=A(\al,\be,\gam)$ is defined
as follows, see \cite[Proposition 2.2]{BenkartRoby}
$$d.v_0=0$$
$$d.v_n=\lambda_{n-1}v_{n-1}, \mbox{for all}\; n\geq 1$$
$$u.v=v_{n+1}.$$

Let ${\frak{h}}=Kud\oplus Kdu$.
We say that an $A$-module $V$ is a weight module if $V=\sum_{\nu\in{\frak{h}}^*}V_{\nu}$,
where $V_{\nu}=\{v\in V|h.v=\nu(h)v\;\mbox{for all}\; h\in{\frak{h}}\}$, and the sum
is over elements in the dual space ${\frak{h}}^*$ of ${\frak{h}}$. If $V_{\nu}\neq 0$, then
$\nu$ is a weight and $V_{\nu}$ is the corresponding weight space. Each weight $\nu$ is
determined by a pair of elements  $(\nu',\nu'')$ of $K$ where
$\nu'=\nu(du)$ and $\nu''=\nu(ud)$,
and we will identify $\nu$ with $(\nu',\nu'')$.

An $A(\al,\be,\gam)$-module $V$ is said to be a highest weight module of weight $\lambda$
if $V$ has a vector $v$ such that $d.v=0$, $du.v=\lambda v$ and $V=A(\al,\be,\gam)v$.
The vector $v$ is said to be a highest weight vector of $V$.

\subsection{\label{novo}} The next Lemma is easily proved by induction,
taking into account the recursive
construction of the $\lambda_i$.

\begin{lem*}
For all $n\in{\Bbb N}$,
$\sigma^{-1}(x-\lambda_n,y-\lambda_{n+1})=(x-\lambda_{n+1},y-\lambda_{n+2})$.
\end{lem*}

We remark that, if $\{v_i|n=0,1,2,...\}$ is the basis of the Verma module, $V(\lambda)$,
$Kv_i\iso R/\sigma^{-i}((x,y-\lambda))$
and we can write the Verma modules using the notation of \cite{Jordan93},
as
$$V(\lambda)\iso\oplus_{i\geq 0}R/\sigma^{-i}((x,y-\lambda)).$$

Also if $dim (L(\lambda))=n$ and we set for a given maximal ideal
${\underline m}$ of $R=K[x,y]$, $L[{\underline m}]=\{a\in
L(\lambda)| {\underline m}a=0\}$, then
$$L(\lambda)=\oplus_{i=0}^{n-1}L[\sigma^{-i}((x,y-\lambda))].$$

If $L[{\underline m}]\neq 0$, we say that ${\underline m}$ is a
weight of $L(\lambda)$.

\subsection{\label{S3}}{\bf Proposition.\;}\;{\em
Every finite dimensional simple $A$-module $M$ is isomorphic to one of the following
\begin{enumerate}
\item[1)] a simple homomorphic image of a Verma module;
\item[2)] 
          $M_P$ where $P$ is a maximal ideal of $R$ which
          has a finite $\sigma$-orbit.
\end{enumerate}}
\begin{proof}
Let $M$ be any finite dimensional simple $A$-module.
The $A_d$-module $A_d\otimes_AM$ is either simple or zero, by
\cite[Theorem 9.17]{GW}.

If $A_d\otimes_A M=0$, $M$ is $d$-torsion. Let $M_0=\{m\in
M|dm=0\}$. Then $(ud)M_0=0$ and as $d^2um=\al dudm+\be ud^2m+\gam
dm=0$, for all $m\in M_0$, we have $(du)M_0\subseteq M_0$. Hence
there is $m\in M_0$ such that $dum=\lambda m$ for $\lambda\in K$.
Now $Am$ is a highest weight module of weight $\lambda$, hence a
homomorphic image of the Verma module $V(\lambda)$,
\cite[Proposition 2.8]{BenkartRoby}. As $M$ is a simple module we
have $M=Am$.

If $A_d\otimes_A M\neq 0$, then as $M$ is simple, $M$ is torsion free and we can identify
$M$ with $A_d\otimes_A M$. The result follows now by Lemma \ref{S2}.
\end{proof}

\begin{note*}
The simple homomorphic images of Verma modules are described in
\cite[Corollary 2.28]{BenkartRoby}.
\end{note*}

\begin{cor*}
Let $A(\al,\be,\gam)$ be a down-up algebra with the parameters
$\al,\be,\gam$ satisfying case (2) of \S \ref{S1} and assume that
$\gam\neq 0$. Then any Verma module $V(\lambda)$ has a unique
maximal submodule $M(\lambda)$. Also any  finite dimensional
simple $A$-module is d-torsion and isomorphic to
$L(\lambda)=V(\lambda)/M(\lambda)$ for some $\lambda$.
\end{cor*}
\begin{proof} Necessary and sufficient conditions for the weight
spaces of $V(\lambda)$ to be one dimensional are given in
\cite[Theorem 2.13]{BenkartRoby}, and in particular these
conditions hold in case (2) of \S \ref{S1} when $\gam\neq 0$. The
statement about Verma modules now follows from \cite[Proposition
2.23]{BenkartRoby}. By Lemma \ref{S3'} any finite dimensional
$A$-module is $d$-torsion, and so has the form $L(\lambda)$ by the
proof of Proposition \ref{S3}.
\end{proof}

\subsection{\label{S4}}
As mentioned in the introduction the algebras $A_{\eta}$ are exactly the down-up
algebras $A(\al,\be,\gam)$ with $\al+\be=1$ and $\gam\neq 0$.
The next result shows that the representation theory of these algebras has certain
similarities with that of $U(osp(1,2))$ and $U(sl_2)$. We return to this topic in section 5.

Note that the algebra $U(sl_2)$ is the only down-up algebra with $\gam\neq 0$ whose
parameters satisfy case (4). We ignore this case below.

The recurrence relation for the $\lambda_n$ is solved explicitly
in \cite[Proposition 2.12]{BenkartRoby}. We use this result below.

\begin{lem*}
Let $A=A_{\eta}$ with $\eta\neq 1$. Then
\begin{enumerate}
 \item[i)] if $\lambda_{n-1}=0$ then $\eta^n\neq 1$;
 \item[ii)] $\lambda_{n-1}=0$ if and only if
 $$\lambda(\eta-1)=-\gam(1-n(\sum_{i=0}^{n-1}\eta^i)^{-1}).$$
\end{enumerate}
\end{lem*}
\begin{proof}
By \cite[Proposition 2.12 (i)]{BenkartRoby},
$\lambda_n=c_1r_1^n+c_2r_2^n+x_n$ where $r_1=1$, $r_2=\eta$,
$x_n=\gam n(1-\eta)^{-1}$, $c_1=(\eta-1)^{-1}[-\lambda -\gamma
+\gam (1-\eta )^{-1}]$ and $c_2=(\eta-1)^{-1}[\eta\lambda +\gamma
-\gam (1-\eta )^{-1}]$. It follows that $\lambda_{n-1}=0$ if and
only if
$$\lambda(\eta^n-1)+\gam(\eta^{n-1}-n)+\gam(1-\eta)^{-1}(1-\eta^{n-1})=0.$$
If $\eta^n=1$, then since $\gam\neq 0$, it follows that
$$\eta^{-1}-n+(1-\eta^{-1})(1-\eta)^{-1}=0$$ or equivalently
$$n(1-\eta)=0$$ so $\eta=1$, a contradiction. This proves $i)$ and
$ii)$ follows by multiplying by
\linebreak
$(\eta-1)/(\eta^n-1)$.
\end{proof}

\begin{prop*}
Let $A=A_{\eta}$ with $\eta\neq 1$.
The only finite dimensional simple modules of dimension
$n\in {\Bbb N}$ are d-torsion modules of the form
$L(\lambda)$ where $n$ is the least
positive integer such that $\lambda$ satisfies $\lambda_{n-1}=0$.
\end{prop*}
\begin{proof}
Let $A$ be a down-up algebra as in the statement of the
proposition. By Corollary \ref{S3} we have that the only finite
dimensional simple $A$-modules are d-torsion and of the form
$L(\lambda)=V(\lambda)/M(\lambda)$.


By construction, the dimension of $L(\lambda)$ is $n$ if and only if $n$ is minimal
with $\lambda_{n-1}=0$, and the result follows.
\end{proof}

It is well known that all simple modules over $U(osp(1,2))$ have odd dimension.
The next result gives a generalization of this fact.

\begin{cor*}
Assume $char(K)=0$.
Let $A=A_{\eta}$ with $\eta\neq 1$ and $\eta$ a primitive $N^{th}$ root of unity. Then
\begin{enumerate}
\item[i)] if $n$ is a multiple of $N$ there are no finite dimensional simple modules of
dimension $n$.
\item[ii)] if $n$ is not a multiple of $N$ there is a unique finite dimensional simple
module of dimension $n$.
\end{enumerate}
\end{cor*}
\begin{proof}
Immediate from the Proposition and the Lemma.
%
\end{proof}

\subsection{\label{dual}} By \cite[\S2]{BenkartRoby}, there is an
antiautomorphism $\tau$ of $A$ given by $\tau(u)=d$, $\tau(d)=u$.
Let ${\cal C}$ be the category of finite dimensional left
$A$-modules. If $M\in {\cal C}$, let $M^*$ be the dual vector
space to $M$. Then $M^*$ is a left $A$-module via
$$(af)(m)=f(\tau(a)m)$$ for $a\in A, f\in M^*, m\in M$. We can now
define a contragradient functor  $(\;)^*$ on the category  ${\cal
C}$ as follows. If $M\in Obj({\cal C})$, let $M^*$ be the dual
vector space to $M$ and if $\psi :M\longrightarrow N$ is a map of
left $A$-modules, then $\psi^*:N^*\longrightarrow M^*$ is the
transpose of $\psi$.

Next suppose $L(\lambda)$ is a finite dimensional highest weight
module with dimension $n$. Then $L(\lambda)$ has a basis
$v_0,\ldots,v_{n-1}$ such that $$uv_i=v_{i+1}\qquad 0\leq i\leq
n-2, uv_{n-1}=0,$$ $$dv_i=\lambda_{i-1}v_{i-1}\qquad 1\leq i\leq
n-1, dv_0=0.$$

Let $f_0,\ldots,f_{n-1}$ be the basis of $L(\lambda)^*$ such that
$f_i(v_j)=\delta_{ij}$ for all $i,j$. Then, $$uf_i=\lambda_i
f_{i+1}\qquad 0\leq i\leq n-2,\qquad uf_{n-1}=0,$$
$$df_i=f_{i-1}\qquad 1\leq i\leq n-1,\qquad df_0=0.$$

In particular $f_0$ is a highest weight vector with weight
$\lambda$ and $L(\lambda)^*\iso L(\lambda)$.

Now suppose that $L(\lambda), L(\mu)$ are finite dimensional
highest weight modules and that $Ext(L(\mu),L(\lambda))\neq 0$.
Then there is a nonsplit exact sequence $$0\longrightarrow
L(\lambda) \longrightarrow M\longrightarrow L(\mu)\longrightarrow
0$$ of $A$-modules. Dualizing we obain $Ext(L(\lambda),L(\mu))\neq
0$.

\begin{cor*}
Suppose that $L(\lambda), L(\mu)$ are finite
dimensional. Then $Ext(L(\lambda),L(\mu))=0$ if and only if
$Ext(L(\mu),L(\lambda))=0$.
\end{cor*}


\section{Down-up modules}

\subsection{\label{Ian1}}
In this subsection and 3.2, $\be$ is allowed to be zero.

Let $A=A(\al,\be,\gam)$ be a down-up algebra and $M$ a left $A$-module.
We define two filtrations on $M$ and view $d$ and $u$ as operators which
move down and up these filtrations (whence the title of this section). The
filtrations need not to be exhaustive.
For any $r$ and $s$ in ${\Bbb N}$,
we define
$$M_r=\{m\in M| d^{r+1}m=0\},$$
$$M^s=\{m\in M| u^{s+1}m=0\},$$
and $$M_r^s=M_r\cap M^s.$$

It is obvious that
$$M_0\subseteq M_1 \subseteq M_2\subseteq\ldots$$
and that $M_r=\{m\in M|dm\in M_{r-1}\}$.

If $\be\neq 0$ then $\cup M_r$ is the $d$-torsion submodule and $\cup M^s$ is the
$u$-torsion submodule of $M$.

\begin{lem*}
For any $r\in{\Bbb N}$, $uM_r\subseteq M_{r+1}$ and $dM^s\subseteq M^{s+1}$.
\end{lem*}
\begin{proof}
Let $m\in M_0$. By $(R1)$ we have $d^2um=0$ and so $uM_0\subseteq
M_1$. Suppose $uM_n\subseteq M_{n+1}$ for all $n< r$. Then
$duM_n\subseteq M_n$ for all $n< r$. If $m\in M_r$, then since
$dudm\in (du)M_{r-1}\subseteq M_{r-1}$, $ud^2m\in
uM_{r-2}\subseteq M_{r-1}$ and $dm\in M_{r-1}$ we have
$$d^2um=(\al dud+\be ud^2+\gam d)m\in M_{r-1}.$$
It follows that
$uM_r\subseteq M_{r+1}$.

The other inclusion is proved in a similar way.
\end{proof}

\begin{cor*}
For any $r,s\in {\Bbb N}$, $dM^s_r\subseteq M^{s+1}_{r-1}$ and
$uM^s_r\subseteq M^{s-1}_{r+1}$.
\end{cor*}
\begin{proof}
Follows easily from Lemma \ref{Ian1}.
\end{proof}

\subsection{\label{Ian4}}
Let $M^{\infty}=\cup M^s$, $M_{\infty}=\cup M_r$ and for each $t\in {\Bbb N}$ let
$M(t)=\sum_{r+s=t} M^s_r$.

It follows easily that the sets $M^{\infty}$, $M_{\infty}$ and $M(t)$ are
$A$-submodules of $M$.

\begin{prop*}
If $M$ is a Noetherian $A$-module such that $M=M^{\infty}=M_{\infty}$ then
$M=M(t)$ for some $t\in{\Bbb N}$.
\end{prop*}
\begin{proof}
Consider the chain of $A$-submodule of $M$, $M(1)\subseteq M(2)\subseteq\ldots$. Choose
$t\in{\Bbb N}$ such that $M(t)=M(t+s)$ for all $s\in{\Bbb N}$. If $m\in M$ then there are
$p,q\in{\Bbb N}$ such that $d^{p+1}m=u^{q+1}m=0$ and this implies that $m\in M(t)$.
\end{proof}

\subsection{\label{Ian6}} Given a module $M$ as in Proposition in \ref{Ian4}, we study conditions
for $M$ to be finite dimensional. From
now on we will assume that $\be\neq 0$.

We define two sequences of elements of $R=K[x,y]$ $$x_0=1, y_0=1$$
$$x_n=\sigma(x_{n-1}x) \quad\mbox{and}\quad
y_n=x\sigma^{-1}(y_{n-1})$$ for any $n\in{\Bbb N}$. We claim that
for any $n\in{\Bbb N}$, $d^nu^n=x_n$.

For $n=0$ this is obvious. The induction step follows from
$$d^{n+1}u^{n+1}=dx_nu=z^{-1}x_nxz=\sigma(x_nx)=x_{n+1}.$$
Similarly we have for all $n\in {\Bbb N}$,  $u^nd^n=y_n$

\begin{lem*}
\begin{itemize}
 \item[i)] If $\sigma^n(x)\notin (x)$ for all $n\in {\Bbb N}$, then $x_t\notin (x)$, for all
$t\in {\Bbb N}$;
 \item[ii)] if $\sigma^{-n}(x)\notin (y)$ for all $n\in {\Bbb N}_0$, then $y_t\notin (y)$, for all
$t\in {\Bbb N}_0$.
\end{itemize}
\end{lem*}
\begin{proof}
For each $t\in {\Bbb N}$ write $x_t=\prod_{i=1}^t\sigma^{i}(x)$
and $y_t=\prod_{i=1}^t\sigma^{-(i-1)}(x)$. The result follows
since $(x)$, $(y)$ are prime ideals of $R$.
\end{proof}

\begin{note*}
\em
We note that the conditions $\sigma^n(x)\notin (x)$ for all $n\geq 1$ and
 $\sigma^{-n}(x)\notin (y)$ for all $n\geq 0$ are equivalent. Indeed if there is
$n\in {\Bbb N}$ such that
\linebreak
$\sigma^n(x)\in (x)$ then $\sigma^n(x)=\lambda x$
for some $\lambda\in K^*$.
So $y=\sigma(x)=\lambda\sigma^{1-n}(x)$ and $\sigma^{1-n}(x)\in (y)$. The
converse follows by a  similar argument.
\end{note*}

\subsection{\label{Ian9}} We now state the main result of this section.

\begin{thm*}
Assume that $K$ is an algebraically closed field and that $\be\neq
0$. Let $M$ be a finitely generated left $A$-module and suppose
that $\sigma^n(x)\notin (x)$ for all $n\in{\Bbb N}$. Then
$A_d\otimes_A M=A_u\otimes M=0$ if and only if $M$ has a finite
filtration whose factor modules are finite dimensional highest
weight modules.

Conversely if every finitely generated left $A$-module which is
both $d$-torsion and $u$-torsion has finite dimension then
$\sigma^n(x)\notin (x)$ for all $n\in{\Bbb N}$.
\end{thm*}
\begin{proof}
Obviously if $M$ has a finite filtration whose modules are finite
dimensional highest weight modules then $M$ is $d$-torsion and
$u$-torsion.

Let $M$ be a finitely generated left $A$-module such that
$A_d\otimes_AM=$
\linebreak
$A_u\otimes_AM=0$. Then $M$ is $d$-torsion and $u$-torsion or equivalently,
$M=M^{\infty}=M_{\infty}$.

As $M$ is a finitely generated left $A$-module and $A$ is Noetherian, so is $M$ and
by Lemma \ref{Ian4} it follows that $M=M(t)$ for some $t$.

Assume that $\sigma^n(x)\notin (x)$ for all $n\in{\Bbb N}$. It is enough to show that $M$ contains a finite dimensional submodule $N$ which is
a highest weight module. For then we can set $N=N_1$ and use the same argument to construct
$0\subset N_1\subset N_2\subset\ldots$ provided $(M/N_i)\neq 0$.

Pick $m\in M_0$ such that
$m\neq 0$. Then $dm=0$. Since $M=M(t)$ we have $u^{t+1}m=0$. By Lemma \ref{Ian4}, we have
$$x_{t+1}m=0.$$
By Lemma \ref{Ian6}, $x_{t+1}\notin (x)$ so
${\ann}_A(m)$ contains a nonzero polynomial in $du$. Therefore $du$ has a nonzero
eigenvector in $M_0$ and the result follows.

Assume that every finitely generated $d$-torsion and $u$-torsion
$A$-module has finite dimension. Fix $n\geq 1$ and set
$I=Au^n+Ad$, $J=I\cap K[x,y]$ and $M=A/I$. Since $\{d^m|m\in{\Bbb
N}\}$ and $\{u^m|m\in {\Bbb N}\}$ are Ore sets, $M$  is
$d$-torsion and $u$-torsion. Thus $M$ has finite dimension.

Since $K[x,y]/J$ embeds in $M$, $J$ has finite codimension in $K[x,y]$. Using the graded
ring structure of $A$, \ref{grading}, it is easily seen that $J$ is the ideal of $K[x,y]$ generated
by $x$ and $x_n$. Hence $x_n\notin (x)$ and so $\sigma^n(x)\notin (x)$.
\end{proof}

%
%
%
%

\subsection{\label{ex}}
Next we give an example where the conclusions of Theorem \ref{Ian9} do not hold.

\begin{ex*}
Let $A=A(0,\be,0)$ and consider the $A$-module $M$ with a basis
$\{m_i,n_i|i\in {\Bbb N}\}$ and such that $um_i=n_i$,
$dn_i=m_{i+1}$, $un_i=0$ and $dm_i=0$.

In this case $d^2M=u^2M=0$ so the relations $d^2u=\be ud^2$, $du^2=\be u^2d$ are obviously
satisfied, and $M$ is generated as an $A$-module by $m_1$.
\end{ex*}


\section{The Isomorphism problem for Down-up Algebras}

Benkart and Roby divide down-up algebras into four classes such that no two algebras from different
classes are isomorphic (see Proposition \ref{iso2} below). Here we solve the isomorphism problem
for three of the classes and make substantial progress on the fourth.

\subsection{\label{iso1}}
First we note the existence of certain isomorphisms and automorphisms.

\begin{lem*}
\begin{itemize}
 \item[i)] If $\be\neq 0$, then $A(\al,\be,\gam)\iso A(-\al\be^{-1},\be^{-1},-\gam\be^{-1})$
 via the map interchanging $d$ and $u$.
 \item[ii)] If $\gam\neq 0$ then $A(\al,\be,\gam)\iso A(\al,\be,1)$.
 \item[iii)] If $A=A(\al,\be,\gam)$, $A'=A(\al',\be',\gam')$ and $\Psi:A\longrightarrow A'$
 is an isomorphism with $\Psi(d)=\lambda d'$, $\Psi(u)=\mu u'$ and $\lambda, \mu\in K^*$
 then $\al'=\al$ and $\be'=\be$.
 \item[iv)] If $A=A(\al,\be,\gam)$, $A'=A(\al',\be',\gam')$ and $\Psi:A\longrightarrow A'$
 is an isomorphism with $\Psi(d)=\lambda u'$, $\Psi(u)=\mu d'$ and $\lambda, \mu\in K^*$
 then $\al'=-\al\be^{-1}$ and $\be'=\be^{-1}$.
\end{itemize}
\end{lem*}
\begin{proof} Straightforward
\end{proof}

\subsection{\label{iso2}}
Next we consider commutative homomorphic images of $A=A(\al,\be,\gam)$.
Let $I=\cap\{J|A/J\;\mbox{is commutative}\}$.
Note that $B=A/I$ is commutative since it is a subdirect product of commutative rings.
Thus $B$ is the largest commutative image of $A$ and $Spec(B)$ should
perhaps be thought of as the largest
commutative subscheme of $Spec(A)$.
The algebra $B$ is given by adding the relations $du=ud$ to the
defining relations for $A$.
Note that the closed points of $Spec(B)$ correspond to the one-dimensional $A$-modules, so we
recover \cite[Theorem 6.1]{BenkartRoby}.

Similarly let $A'=A(\al',\be',\gam')$ be another down-up algebra and $I'$ be the unique
smallest ideal of $A'$ such that $A'/I'$ is commutative. Suppose there exists an isomorphism
$\Psi$ from $A$ onto $A'$. It is easily seen that $\Psi(I)=I'$. Moreover
$$\Psi(\sum\{P|P\;\mbox{minimal over}\; I\})=\sum\{P'|P'\;\mbox{minimal over}\; I'\}.$$


The images $a,b$ of $d,u$ in $B$ satisfy
               $$(R3)\qquad a(ab(1-\al-\be)-\gam)=0$$
               $$(R4)\qquad b(ab(1-\al-\be)-\gam)=0.$$
Thus we obtain.

\begin{prop*}
The largest commutative homomorphic image $B$ of $A$ is the factor ring of
the commutative polynomial ring $K[a,b]$ defined by
$(R3),(R4)$. In particular one of the following cases holds
\begin{itemize}
 \item[(a)] $\gam=0$, $\al+\be=1$ and $B=K[a,b]$;
 \item[(b)] $\gam=0$, $\al+\be\neq 1$, $B=K[a,b]/(a^2b,ab^2)$ and the primes of $A$ minimal
 over
 $I$ are $(d)$ and $(u)$;
 \item[(c)] If $\gam \neq 0$, $\al+\be\neq 1$ then the primes of $A$ minimal over $I$
 have the form  $M_0=(d,u)$ and $P={\Ker}(\Phi)$ where $\Phi:A\longrightarrow K[v,v^{-1}]$ is
 defined by $\Phi(u)=(1-\al-\be)^{-1}v$ and $\Phi(d)=\gam v^{-1}$.
    \item[(d)] If $\gam\neq 0$, $\al+\be=1$ then $B=K$ and $I=(d,u)$.
\end{itemize}
\end{prop*}

\subsection{\label{iso3}}
We say that the down-up algebra $A=A(\al,\be,\gam)$ has  type (a) (resp.  (b),
 (c), (d)) if the parameters $\al,\be,\gam$ satisfy condition (a) (resp. (b),
(c), (d)) of Proposition \ref{iso2}.

Let $A=A(\al,\be,\gam)$ and $A'=A(\al',\be',\gam')$ be down-up
algebras of the same type and let $d',u'$ be the generators of
$A'$. If $A$ and $A'$ have type (a) or (d) and $\eta\neq 1\neq
\eta'$, let $w_2$ and $w_2'$ be the elements constructed in case
(2) of \ref{S1}. If $A$ and $A'$ have type (c), let $P$ be the
ideal of $A$ defined in Proposition \ref{iso2} and $P'$ the
corresponding ideal of $A'$.

\begin{cor*}
With the notation as above assume that $A$ and $A'$ are isomorphic via $\Psi$. We have
\begin{itemize}
 \item[i)] If $A$, $A'$ have type (a), then $\Psi(w_2)=(w_2')$;
 \item[ii)] If $A$, $A'$ have type (b), then $\Psi(d,u)=(d',u')$;
 \item[iii)] If $A$, $A'$ have type (c), then $\Psi(P)=P'$.
\end{itemize}
\end{cor*}
\begin{proof} If $A$ and $A'$ have type (a) then using the decomposition in
\ref{grading}, we see that $w_2=-x+y$ generates the ideal $I$.
Hence (i) follows from the remarks before Proposition \ref{iso2}.
The proofs of the remaining statements are similar.
\end{proof}

\subsection{\label{iso3b}}
Before stating our main result on the isomorphism problem, it is
worth commenting on the geometry of one-dimensional
representations, for algebras of type (c). The maximal ideals of
$K[v,v^{-1}]$ have the form $(v-\mu)$, $\mu\in K^*$, so we have
homomorphisms $\Phi_{\mu}:A\longrightarrow K$ given by
$$\Phi_{\mu}(u)=(1-\al-\be){\mu},\; \Phi_{\mu}(d)=\gam\mu^{-1}.$$
Set $M_{\mu}=Ker(\Phi_{\mu})$. Then as in \cite{BenkartRoby}, the
one dimensional modules are indexed by $K^*\cup\{0\}=K$. However
from Proposition \ref{iso2} we might expect the ideal $M_0$ to
behave differently from the other $M_{\mu}$. Indeed we have

\begin{lem*}
For algebras of type (c), $M^2_{\mu}=M_{\mu}$ if and only if $\mu=0$.
\end{lem*}
\begin{proof}
Since $M_0=(d,u)$ and $\gam\neq 0$, relations $(R1), (R2)$ imply that $M_0^3=M_0$. On the
other hand if $\mu\neq 0$ then $(v-\mu)\neq (v-\mu)^2$ in the commutative ring $K[v,v^{-1}]$,
so it follows that $M_{\mu}\neq M_{\mu}^2$.
\end{proof}

\subsection{\label{iso4}}
It follows from \cite[Corollary 6.2]{BenkartRoby} that if two down-up algebras
are isomorphic then they have the same type.
The next result gives a partial answer to
\cite[Problem (h)]{BenkartRoby} and a partial converse of Lemma \ref{iso1}.

\begin{thm*}
Suppose that $A=A(\al,\be,\gam)$ and $A'=A(\al',\be',\gam')$ are Noetherian down-up algebras
of the same type. If both have type (d) assume also that $char(K)=0$.
Then $A\iso A'$ if and only if

\begin{tabular}{rll}
 & 1) & $\gam=0$ if and only if $\gam'=0$ and\\
either & 2) & $\al=\al'$, $\be=\be'$\\
or & 3) & $\al=-\al\be^{-1}$, $\be'=\be^{-1}$.\\
\end{tabular}
\end{thm*}

\subsection{\label{iso5}} Assume that $A=A(\al,\be,\gam)$ with $\al+\be =1$,
and $\eta\neq 1$. Then case 2 of \ref{S1} holds and we set
$w=w_2$, so that $\sigma(w)=\eta w$.

\begin{lem*}
 The set $\{a\in A|aw=\eta^mwa\}$ equals\\
 \begin{tabular}{ccl}
 A(m) & if & $\eta$ is not a root of unity\\
 $\bigoplus \{A(m')|m'\equiv m\; (mod\, n)\}$ & if & $\eta$ is a primitive $n^{th}$
 root of unity.
 \end{tabular}
\end{lem*}
\begin{proof}
This is
proved by computation using the decomposition in \ref{grading}.
\end{proof}

\subsection{\label{iso7}} Now assume that $A$ is a down-up algebra of type (a) or
(b).

\begin{lem*}
Let $A$ and $A'$ be down-up algebras both of type (a) with
$\eta\neq 1\neq \eta'$ or of type (b). Assume that $A$ and $A'$
are isomorphic via $\Psi$. Then $\Psi(d,u)\subseteq (d',u')$.
\end{lem*}
\begin{proof}
In type (b) this follows directly from Corollary \ref{iso3}.
Suppose that $A, A'$ have type (a) and that $\eta\neq 1\neq\eta'$.
By Corollary \ref{iso3} and the fact that $A$ has only trivial
units we have $\Psi(w)=\lambda w'$ for some $\lambda\in K^*$,
where $w'\in A'$ is defined in a similar manner to $w$. Applying
$\Psi$ to the equation $dw=\eta wd$, we get by Lemma \ref{iso5},
$\Psi(d)\in \sum \{A'(m)|(\eta')^m=\eta\}\subseteq (d',u')$.
Similarly $\Psi(u)\in (d',u')$.
\end{proof}

\subsection{\label{iso8}} {\bf Proof of Theorem \ref{iso4} for type (a) and type (b)}.

Adopting some terminology from group theory we say that a subset $X$ of $A=A(\al,\be,\gam)$ is
{\em characteristic} if $\Psi(X)=X$ for all $\Psi\in Aut(A)$.

Assume that $A$ and $A'$ are down-up algebras both of type (a) or
type (b). Assume as well that $\be,\be'\neq -1$ if $A$ and $A'$
are of type (a). Since $\gam=0$, $A=\oplus_{m\in {\Bbb N}_0}A[m]$
is a graded algebra with $A[0]=K$ and $A[1]=span\{d,u\}$ and
similarly for $A'$. Set $A_n=\oplus_{m\geq n}A[m]$. Then
$(d,u)=A_1$ and $(d,u)^n=A_1^n=A_n$. Hence if
$\Psi:A\longrightarrow A'$ is an isomorphism we have
$\Psi(A_n)\subseteq A'_n$ by Lemma \ref{iso7}. Thus $\Psi$ is an
isomorphism of graded algebras. As noted in \cite[Theorem 4.1 and
Lemma 4.2]{KMP} $A$ is Auslander regular of global dimension 3.
Since $A$ has GK-dimension 3 by \cite[4.2]{Benkart} it follows
from \cite[Theorem 6.3]{L} that $A$ is Artin-Schelter regular.
Therefore the isomorphism type of $A$ as a graded algebra is
determined in \cite{AS}. However instead of appealing to \cite{AS}
we can now complete the proof with a short calculation. Since the
relations are of degree three we work $mod A'_4$. Note that the
algebra $A/A_4$ inherits the ${\Bbb Z}$-grading of \ref{grading}.
Write $$\Psi(d)=bd'+cu'$$ $$\Psi(u)=rd'+su'$$ $mod A'_2$.
Obviously $\Delta = bs-cr\neq 0$. Suppose first that $A$, $A'$
have type (b), so $\al+\be\neq 1$. Applying $\Psi$ to the relation
$d^2u-\al dud - \be ud^2=0$ and looking at the terms of degree $3$
and $-3$ we see that $br=cs=0$. Since $\Delta\neq 0$ this gives
two possibilities. Either $b\neq 0\neq s$ and $c=r=0$ or
$c\neq0\neq r$ and $b=s=0$. It is easily seen that we have one of
the two statements of the Theorem.

Now suppose that $A$ and $A'$ have type (a) and $\be,\be'\neq -1$.
Applying $\Psi$ to the relation
$d^2u-\al dud - \be ud^2$ and cancelling $\Delta$ we obtain
$$0=b[(d')^2u'-\al d'u'd'-\be u'(d')^2]+c[\be d'(u')^2+\al u'd'u'- (u')^2d']$$
$mod A'_4$.
Comparing to the relations in $A'$ gives the result.

Finally suppose that $A,A'$ are down up algebras of type (a),
$\eta =1$ and $\be'\neq -1$. By Lemma \ref{iso7}, $A'$ has a
characteristic ideal of the form $(d',u')$. Since $\eta=1$,
$A=U({\frak{h}})$, the enveloping algebra of the Heisenberg
algebra ${\frak h}$. Now ${\frak h}$ has basis $\{x,y,z\}$ such
that $[x,y]=z$ is central in ${\frak h}$. By \cite[Theorem
6.1]{BenkartRoby} any codimension one ideal in $A$ has the form
$(x-a,y-b,z)$. It is easy to see  that $Aut(A)$ acts transitively
on maximal ideals of codimension 1, hence there is no
characteristic ideals of codimension 1. This proves Theorem
\ref{iso4} for type (a) and (b).

\subsection{\label{iso10}} {\bf Proof of Theorem \ref{iso4} for type (c)}.

To solve the isomorphism problem for down-up algebras of type (c)
we consider the class of bimodules over $C=K[v,v^{-1}]$ which are
free of rank $n$ on the left. Let $F$ be such a bimodule with
basis $e_1,\ldots,e_n$. We can define a ${\Bbb Z}$-grading
$\{F(m)|m\in {\Bbb Z}\}$ on $F$ by $F(m)=\sum_{i=1}^nKv^me_i$. We
assume that the right $C$-action preserves this grading, that is
$e_iv\in F(1)$ for all $i$. Then we can write
$$e_iv=\sum_{j}p_{ij}ve_j$$ for some $n\times n$ matrix
$P=(p_{ij})$ with $p_{ij}\in K$. If $\phi$ is an automorphism of
$F$ and $$\phi(e_i)=\sum q_{ij}e_j=f_i$$ with $q_{ij}\in C$, then
the right action of $v$ on the basis $f_1,\ldots f_n$ is
determined by the matrix $QPQ^{-1}$. Thus the isomorphism class of
the bimodule $F$ is determined by the conjugacy class of $P$ under
the action of $GL_n(C)$.

\begin{lem*}
Suppose $A=A(\al,\be,\gam)$ is a down-up algebra of type (c), let
\linebreak $\Phi:A\longrightarrow C=K[v,v^{-1}]$ be the
homomorphism described in Proposition \ref{iso2} and
$P={\Ker}(\Phi)$. Then $P/P^2$ is free of rank 2 as a left and
right $C$-module. In addition one of the following holds:
\begin{itemize}
 \item[i)] $\al^2+4\be\neq 0$. Let $\lambda_1, \lambda_2$ be the roots of the polynomial
 $f(\lambda)$ and
 $$w_i=\be(\lambda_i-1)ud-\lambda_i(\lambda_i-1)du+\gam\lambda_{i}.$$
 Then $w_1,w_2$ are normal elements of $A$ and $P=(w_1,w_2)$. The bimodule $P/P^2$ is
 free of rank two on the left and right as a $C$-module whose isomorphism class is
 determined by the conjugacy class of the matrix
 $$ \left[\begin{array}{rr}
          \lambda_1 & 0\\
           0        & \lambda_2
          \end{array}
    \right]$$
 \item[ii)] $\al^2+4\be =0$. Then $\al/2$ is a double root of $f(\lambda)$. Let
            $$w_1=(2\be+\al)ud+(\al-2)du+2\gam$$
            $$w_2=2(du-ud).$$
            Then $w_1$ is normal in $A$, the image of $w_2$ is normal in $A/(w_1)$ and
            \linebreak
            $P=(w_1,w_2)$.
            The bimodule $P/P^2$ is free of rank two on the left and right as a
            $C$-module whose isomorphism class is
            determined by the conjugacy class of the matrix
            $$ \left[\begin{array}{rr}
                     \al/2  & 1\\
                      0     & \al/2
               \end{array}
               \right]$$
\end{itemize}
\end{lem*}
\begin{proof}
We prove only part ii). The proof of part i) is similar. The fact
that $w_1$ is normal in $A$ and $w_2$ normal $mod (w_1)$ follows
from \ref{S1}. A short computation shows that $w_1, w_2\in
{\Ker}(\Phi)$. Using the decomposition of $A$ as a ${\Bbb
Z}$-graded ring in \ref{grading} we see that $P/P^2$ is free as a
left and right $C$-module with basis ${\overline w}_i=w_i+P^2$,
$i=1,2$. From  \ref{S1}, we obtain $${\overline
w}_1v=v(\al/2{\overline w}_1+{\overline w}_2),$$ $${\overline
w}_2v=(\al/2)v{\overline w}_2$$ and the result follows.
\end{proof}

Before concluding the proof of Theorem \ref{iso4} in this case we need another definition.
Suppose that $M$ is a $C$-bimodule and $\mu\in Aut(C)$. The $C$-bimodule twisted by $\mu$
has the same underlying vector space as $M$, and has bimodule structure maps
$C\times M\rightarrow M$, $M\times C\rightarrow M$ given by
$(c,m)\mapsto \mu(c)m$ and $(m,c)\mapsto m\mu(c).$

Now suppose $\Psi:A\rightarrow A'$ is an isomorphism of down-up
algebras of type (c) and let $\Phi:A\rightarrow C$,
$\Phi':A'\rightarrow C$ be the maps described in Proposition
\ref{iso2} (c). Let $P={\Ker}(\Phi)$, $P'={\Ker}(\Phi')$ and write
$\phi:A/P\rightarrow C$, $\phi':A'/P'\rightarrow C$ for the
induced isomorphisms. There is an automorphism $\mu$ of $C$
satisfying $\mu\Phi=\Phi'\Psi$. Since $\Psi(P)=P'$, by \ref{iso3},
$\Psi$ induces a linear isomorphism from $P/P^2$ to $P'/(P')^2$.
If we regard $P/P^2$ as a $C$-bimodule via $\phi^{-1}$ and
$P'/(P')^2$ as a $C$-bimodule via $(\phi')^{-1}\mu$ then the above
map is an isomorphism of $C$-bimodules.

Now the isomorphism type of $P/P^2$ (resp. $P'/(P')^2$) is determined by the conjugacy classes of a matrix
$J$, respectively $J'$, as in the Lemma. It follows that these bimodules
structures can be obtained from one another twisting by $\mu,\mu^{-1}\in Aut(C)$.

There are now two possibilities. If $\mu(v)=\lambda v$ for some $\lambda\in K^*$, then
$J,J'$ are conjugate and we have conclusion (2) of Theorem \ref{iso4} while if
$\mu(v)=\lambda v^{-1}$ for some $\lambda\in K^*$, then $J^{-1}$ and
$J'$ are conjugate and we have conclusion (3).

\subsection{\label{iso3a}}
We want to establish an analogue of Corollary \ref{iso3} (a) for down-up algebras
of type (d).

Let $A$ be a down-up algebra of type (d) and assume $\eta\neq 1$.
Let \linebreak $w_2=du-ud-\gam (\eta-1)^{-1}$. From section
\ref{S1} we have $dw_2=\eta w_2 d$ and $w_2 u=\eta uw_2$. In
particular $w_2$ is a normal element. Clearly $A/(w_2)$ is
isomorphic to the first Weyl algebra since $\gam \neq 0$. Thus
$(w_2)$ is a completely prime ideal of $A$.

\begin{lem*}
Let $A$ be a down-up algebra of type (d) and assume that $char(K)=0$ and $\eta\neq 1$.
If $P$ is a completely prime ideal of $A$ such that $A/P$ has infinite dimension over $K$,
then $P=(w_2)$.
\end{lem*}
\begin{proof}
Suppose $P$ is a completely prime ideal such that $A/P$ has infinite dimension over $K$.
Then $(d,u)\nsubseteq P$ so assume that $d\notin P$ and localize at $d$.
If $u\notin P$ we localize at $u$ instead and use a similar argument.
Then $Q=P_d$ is a nontrivial ideal of $S=A_d=R[z,z^{-1};\sigma]$.
If $Q\cap R=0$ we can localize at $R\backslash\{0\}$ to obtain a nontrivial ideal in
$F[z,z^{-1};\sigma]$ where $F=Fract(R)$. However it follows from section \ref{S1}
and the assumption that $char(F)=0$, that $\sigma$ has infinite order.
Hence by \cite[Theorem 1.8.5]{M&R}, $F[z,z^{-1};\sigma]$ is a simple ring.

This contradiction shows that $I=Q\cap R\neq 0$. Now by \ref{S1}, $R=K[w_1,w_2]$ where
$\sigma (w_1)=w_1+\gam$ and $\sigma (w_2)=-\be w_2$.

Now $R/I$ embeds in $S/Q$ which is a domain, so $I$ is prime and clearly $\sigma$-invariant.
By choosing a polynomial of least degree in $w_2$ with coefficients in $K[w_1]$ we see
that $(w_2)\subseteq I$ and the lemma follows easily.
\end{proof}

\begin{cor*}
Assume that $char(K)=0$ and $\eta\neq 1\neq\eta'$. Let
$A=A(\al,\be,\gam)$ and $A'=A(\al',\be',\gam')$ be isomorphic
down-up algebras of type (d) via the isomorphism $\Psi$. Let $w_2$
be the element of $A$ defined above and $w_2'$ the corresponding
element of $A'$. Then $\Psi(w_2)=(w_2')$.
\end{cor*}

\subsection{\label{iso11}} {\bf Proof of Theorem \ref{iso4} for type (d) and char(K)=0}.

Suppose first that $\al^2+4\be\neq 0$ and $(\al')^2+4\be'\neq 0$. Using the notation
of Corollary \ref{iso3a} we have
$$\Psi(w_2)=\lambda w'_2$$
for some $\lambda\in K^*$ since $A'$ has no nontrivial units.
By Corollary \ref{S4}, $\eta$ is an $n^{th}$ root of unity if and only if
$\eta'$ is an $n^{th}$ root of unity.

If $\eta$ and $\eta'$ are not roots of unity, applying $\Psi$ to
the equation $dw_2=\eta w_2d$ and using Lemma \ref{iso5} gives
$\Psi(d)\in\sum\{A'(m)|m\in{\Bbb Z}, (\eta')^m=\eta\}$. There can
be at most one $m\in{\Bbb Z}$ such that $(\eta')^m=\eta$. If $m>
0$ then $\Psi(d)\in (d')^mR\subseteq (d')$ and $\Psi(u)\in
(u')^mR\subseteq (u')$. If $m< 0$, $\Psi(d)\in (u')^mR\subseteq
(u')$ and $\Psi(u)\in (d')^mR\subseteq (d')$. As $A$ has no
nontrivial units, the result follows from Lemma \ref{iso1}.

Assume now that $\eta$, $\eta'$ are primitive $n^{th}$ roots of
unity, $n>1$ and that $\gam =1$. If $n=2$, then $\eta=\eta'=-1$,
so we can assume $n>2$. By Proposition \ref{S4}, there are unique
simple modules $L_1$, $L_2$ of dimension $1, 2$ and these modules
are weight modules with highest weights  $0$ and
$-\al^{-1}=(\be-1)^{-1}$, respectively. Since
$w_2=du-ud-(\be+1)^{-1}$ we see that $w_2$ acts on the highest
weight vectors of $L_1$, $L_2$ by the eigenvalues $-(\be+1)^{-1}$
and $2(\be^2-1)^{-1}$, respectively. Since $w_2^n$ is central it
acts as the $n^{th}$ power of these eigenvalues on $L_1, L_2$.

Now suppose $\Psi:A=A(\al,\be,1)\longrightarrow A'=A(\al',\be',1)$ is an isomorphism and
define $w_2'\in A'$ analogously to the way $w_2\in A$ is defined. By Corollary \ref{iso3a}
$\Psi (w_2)=\lambda w_2'$. Thus comparing the actions of $w_2^n$ and $(w_2')^n$ we get
$$\lambda^n=(\frac{(\be')^2-1}{\be^2-1})^n=(\frac{\be' +1}{\be +1})^n.$$
Hence $(\be-1)^n=(\be'-1)^n$. This implies that $|\be-1|=|\be'-1|$ and then by writing
$\be=e^{i\th}$ that $\be'=\be$ or $\be'=\be^{-1}$.

It remains to rule out the possibility that $A$ and $A'$ are
isomorphic down-up algebras of type (d) and $\al^2+4\be=0$,
$(\al')^2+4\be'\neq 0$. Hence $A\iso A(2,-1,1)\iso U(sl_2))$. If
$\eta'$ is a root of unity we obtain a contradiction by Corollary
\ref{S4}. Thus $\eta'$ is not a root of unity and by \cite[Theorem
4.0.2]{K} or \cite[Theorem 1.3]{Zhao}, the centre of $A'$,
$Z(A')=K\neq Z(A)$, so $A\ncong A'$.

\section{semisimplicity}

In this section we assume $char(K)=0$ and we study down-up
algebras of type (d). For brevity we use the notation $A_{\eta}$
to refer to these algebras, since as observed in the introduction
the $A_{\eta}$ are exactly the down-up algebras of type (d).
However we continue to use the generators $d,u$ and relations
$(R1)$, $(R2)$ rather than the generators $h,e,f$ to preserve
continuity. Recall that $\al=\eta +1$, $\be=-\eta$ and $\gam\neq
0$. Our main result is that all finitely generated $A_{\eta}$
modules are semisimple if $\eta$ is a root of unity.

Given any $\lambda$, we define the following maximal ideals of $R=K[x,y]$,
$$J_{\lambda}=(x,y-\lambda),$$
$$H_{\lambda}=(x-(\lambda-\gam)\be^{-1},y).$$

It is trivial to confirm that $H_{\lambda}=\sigma(J_{\lambda})$ and that
$dJ_{\lambda}=H_{\lambda}d$.

\subsection{\label{semisimple}}
The proof of the next result is adapted from \cite[Theorem 5.2]{Jordan93}.

\begin{prop*}
Let $A=A_{\eta}$, $\lambda\neq\mu$ and $$0\longrightarrow
L(\lambda)\longrightarrow M \longrightarrow L(\mu)\longrightarrow
0$$ a nonsplit short exact sequence of finite dimensional
$A$-modules. Then one of the following occurs
\begin{itemize}
 \item[i)] $M$ is an epimorphic image of $V(\mu)$;
 \item[ii)] there is a nonzero $w\in L(\lambda)$ such that $H_{\mu}w=0$.
\end{itemize}
\end{prop*}
\begin{proof}
As in \cite[Theorem 5.2]{Jordan93} we can choose $v\in M\backslash
L(\lambda)$ such that the image in $L(\mu)$ is a highest weight
vector in $L(\mu)$ and that $J_{\mu}^2v=0$. Now $J_{\mu}v\subseteq
L(\lambda)$ and either $dJ_{\mu}v=0$ or $dJ_{\mu}v\neq 0$.

Assume that $dJ_{\mu}v\neq 0$. Then since $J^2_{\mu}v=0$ we get
$0=dJ^2_{\mu}v=H_{\mu}(dJ_{\mu}v)=0$ so $ii)$ holds.

If $dJ_{\mu}v=0$, then $J_{\mu}v=0$, since otherwise we would find
$w\in J_{\mu}v$, $w\neq 0$ and then $J_{\mu}w=dw=0$ would
contradict $\lambda\neq \mu$. If also $dv=0$ then $i)$ holds. If
$J_{\mu}v=0\neq dv$, then since $H_{\mu}dv=dJ_{\mu}v=0$, $ii)$
holds.
\end{proof}

\subsection{\label{weights2}} Easy calculations show that there is
no nonzero integer $k$ such that
$\sigma^k(J_{\lambda})=J_{\lambda}$.

\begin{lem*}
If there is $k\in {\Bbb Z}$ such that
$\sigma^k(J_{\lambda})=J_{\lambda}$, then $k=0$.
\end{lem*}
\begin{proof}
Let $k\in {\Bbb Z}$ be such that
$$J_{\lambda}=\sigma^{k}(J_{\lambda})=(x-\lambda_{-k-1},y-\lambda_{-k}).$$
Hence $\lambda_{-1}=\lambda_{-k-1}$, $\lambda_0=\lambda_{-k}$, and
by \cite[Theorem 2.23]{BenkartRoby} we conclude that $k=0$.
\end{proof}

\subsection{\label{itself}} We prove that no
nonsplit extensions of $L(\lambda)$ by itself ever occur for $A_{\eta}$.

\begin{lem*}
Let $A=A_{\eta}$ such that $\eta\neq 1$. Assume that
$dim(L(\lambda))=n$. Then $$(x,\sigma^n(x))=(x,y-\lambda).$$
\end{lem*}
\begin{proof}
Let $\{v_i\}$ be a basis for the Verma module $V(\lambda)$. As
$dim(L(\lambda))=n$, $\lambda_i\neq 0$ for all $i<n-1$ and
$\lambda_{n-1}=0$. We have that $J_{\lambda}={\ann}_R(v_0)$ and
$\sigma^{-n}(J_{\lambda})={\ann}(v_n)$. Since
$\sigma^{-n}(J_{\lambda})=(x,y-\lambda_{n})$, $\sigma^n(x)\in
J_{\lambda}$. Hence $(x,\sigma^n(x))\subseteq (x,y-\lambda).$
Combined with the fact that $\sigma$ stabilizes $span\{1,x,y\}$,
this implies that we can write $\sigma^n(x)=ax+b(y-\lambda)$, for
some $a,b\in K$. If $b=0$, then $\sigma^n(x)=ax$. Thus $a$ is
either $1$ or $\eta^n$ and a short calculation shows that $x\in
span\{w_2,1\}$, a contradiction. Hence $y-\lambda\in
(x,\sigma^n(x))$
\end{proof}

The proof of next result is adapted from \cite[Theorem 5.4]{Jordan93}

\begin{prop*}
Let $A=A_{\eta}$ and $\eta\neq 1$. There are no nonsplit short
exact sequences of $A_{\eta}$-modules of the form
$$0\longrightarrow L\longrightarrow M \longrightarrow
N\longrightarrow 0$$ with $L\iso N\iso L(\lambda)$ finite
dimensional.
\end{prop*}
\begin{proof}
Assume that there are sequences as above. First suppose that
$dim(L(\lambda))=n>1$. As in \cite[Theorem 5.2]{Jordan93} choose
$v\in M\backslash L$ such that the image in $N$ is a highest
weight vector and such that $J_{\lambda}^2v=0$. By Lemma
\ref{novo} the weights of $L(\lambda)$ are
$\sigma^{-i}(J_{\lambda})$ with $i=0,\ldots,n-1$. Thus by Lemma
\ref{weights2}, $H_{\lambda}=\sigma(J_{\lambda})$ and
$\sigma^{-n}(J_{\lambda})$ are not weights of $L(\lambda)$. On the
other hand $H_{\lambda}^2dv=dJ_{\lambda}^2v=0$, so $dv=0$.
Similarly since $dim(L({\lambda}))=n$ we have $u^nv\in L$, and
$J_{\lambda}^2v=0$ implies that
$0=u^nJ_{\lambda}^2v=\sigma^{-n}(J_{\lambda}^2)u^nv$. Thus
$u^nv=0$.




Now $u^{n-1}\sigma^n(x)v=\sigma(x)u^{n-1}v=yu^{n-1}v=du^nv=0$ and
also $d\sigma^n(x)v=\sigma^{n+1}(x)dv=0$. Since
$\sigma^n(x)\subseteq J_{\lambda}$, $\sigma^n(x)v\in L$. Since
$n>1$ an easy computation shows that the only element of $L$
annihilated by $u^{n-1}$ and by $d$ is the zero element. Hence
$\sigma^{n}(x)v=0$. Since $(x,\sigma^n(x))=J_{\lambda}$, we have
$J_{\lambda}v=0$. Since $L$ and $M$ are simple, we have that
$M=Av$ and hence $M$ is a highest weight module of weight
$\lambda$. As $u^nv=0$, we conclude that $M\iso L(\lambda)$ and
hence such a nonsplit exact sequence can not occur.
%

Finally suppose that $dim(L(\lambda))=1$, that is $\lambda =0$. Since
$L(0)\iso A/(d,u)$, we have $(d,u)^2M=0$. As in Lemma \ref{iso3}, we have
$(d,u)^2=(d,u)$. Hence $M$ is a module over the field $A/(d,u)$,
so the sequence splits, a contradiction.
\end{proof}

\subsection{\label{final}}If $A$ is a down-up algebra such that
all finite dimensional $A$-modules are semisimple then $A$ has type
$(d)$ since otherwise by Proposition \ref{iso2} $A$ has a commutative
image which is not simple Artinian. Set
$$X_{m,n}=\{\eta\in K^*| \eta^m\neq 1\neq\eta^n\;\mbox{and}\;
n(\eta^m-1)=m(\eta^n-1)\}$$
for $m,n\geq 1$.

\begin{thm*}
Let $A=A_{\eta}$ with $\eta\neq 1$. Then the
following
\begin{enumerate}
 \item[i)] all finite dimensional $A_\eta$-modules are semisimple;
 \item[ii)] every Verma module for $A_\eta$ have composition length
 $\leq 2$;
 \item[iii)] $\eta\notin X_{m,n}$, for all $m\neq n$.
\end{enumerate}
\end{thm*}
\begin{proof}
We first show the equivalence of conditions $i)$ and $ii)$. By
\cite[Theorem 2.13 and Proposition 2.23]{BenkartRoby} any
submodule of $V(\lambda)$ has the form $N=span\{v_j|j\geq n\}$ for
some $n>0$. It follows that $V(\lambda)$ has length $\leq 2$ if
and only if every finite dimensional homomorphic image of
$V(\lambda)$ is simple. Thus $i)$ implies $ii)$. Conversely
suppose that $ii)$ holds and there is a nonsplit exact sequence
$$0\longrightarrow L(\lambda)\longrightarrow M\longrightarrow
L(\mu)\longrightarrow 0$$ with $M$ finite dimensional. By
Proposition \ref{itself} $\lambda\ne \mu$. Since all Verma modules
have composition length $\leq 2$, case (i) of Proposition
\ref{semisimple} cannot occur. Thus $\sigma(J_{\mu})$ is a weight
of $L(\lambda)$, that is $\sigma^i(J_{\mu})=J_{\lambda}$ for some
$i>0$. Similarly by dualizing the above exact sequence and
applying Proposition \ref{semisimple} we get
$\sigma^j(J_{\lambda})=J_{\mu}$ for some $j>0$. Hence
$J_{\lambda}=\sigma^{i+j}(J_{\lambda})$, but this contradicts
Lemma \ref{weights2}.

Next we prove $iii)$ implies $ii)$ by showing that if $V(\lambda)$
is a Verma module over $A_{\eta}$ of length $>2$, then $\eta\in
X_{m,n}$, for some $m\neq n$. By the description of the submodules
of $V(\lambda)$ from \cite{BenkartRoby} cited above we have
$\lambda_{m-1}=\lambda_{n-1}=0$ for some $m>n>0$. Thus by Lemma
\ref{S4}, $\eta^m\neq 1\neq\eta^n$ and $$\lambda(\eta
-1)=-\gam(1-n(\sum_{i=0}^{n-1}\eta^i)^{-1})$$ $$\lambda(\eta
-1)=-\gam(1-m(\sum_{i=0}^{m-1}\eta^i)^{-1}).$$ Forming the
difference between these equations we get
$$\gam[m(\sum_{i=0}^{m-1}\eta^i)^{-1}-n(\sum_{i=0}^{n-1}\eta^i)^{-1}]=0.$$
As $\gam\neq 0$ $$(\eta^n-1)m=(\eta^m-1)n.$$ Thus $\eta\in
X_{m,n}.$

Finally, if $\eta\in X_{m,n}$ and $\lambda$ is chosen so that
$$\lambda(\eta-1)=-\gam(1-m(\sum_{i=0}^{m-1}\eta^i)^{-1})$$ it
follows that $\lambda_{m-1}=\lambda_{n-1}=0$ and the Verma module
$V(\lambda)$ has length $>2$.
\end{proof}

In the preliminary version of this paper, we proved semisimplicity
of finite dimensional $A_{\eta}$-modules when $\eta$ was a root of
unity. The present version of Theorem \ref{final} involves only
minor changes to the earlier proof. We thank D. Jordan for
pointing out the equivalence of conditions $i)$ and $iii)$.

\subsection{\label{end}} Finally we show that

\begin{lem*}
If $\eta$ is a root of unity or $K={\Bbb C}$ and $|\eta|=1$, then
$\eta\in X_{m,n}$ implies $m=n$.
\end{lem*}
\begin{proof}
Suppose $\eta^m\neq 1$, $\eta^n\neq 1$ and
$$(\eta^n-1)m=(\eta^m-1)n\qquad\qquad (*).$$ Again this is an
equation in ${\Bbb Q}(\eta)$ which we can identify with a subfield
of ${\Bbb C}$. Then $(\eta^m-1)/(\eta^n-1)=mn^{-1}\in {\Bbb R}$.
Consideration of the imaginary part of this expression shows that
$$sin(m\th)(cos(n\th)-1)=sin(n\th)(cos(m\th)-1)$$ where
$\eta=e^{i\th}$. If $sin(m\th)=0$, then since $\eta^m\neq1$ we get
$m\th+\pi\in 2\pi{\Bbb Z}$, and $cos(m\th)=-1$. Thus
$sin(n\th)=0$, $\eta^m=\eta^n=-1$ and $(*)$ forces $m=n$. Hence we
can assume that $sin(m\th)\neq 0\neq sin(n\th)$. Let
$g(x)=(cos(x)-1)/sin(x)$, so that $g(m\th)=g(n\th)$. Then $g(x)$
is decreasing on $(-\pi, \pi)$ so $(m-n)\th\in 2\pi{\Bbb Z}$.
Therefore $(*)$ forces $m=n$.
\end{proof}

\begin{prop*}
Any Verma module over $A_{\eta}$ has length $\leq 3$.
\end{prop*}
\begin{proof}
If the result is false then from the description of the submodules
of $V(\lambda)$, we can find positive integers $m<n<p$ such that
$\lambda_{m-1}=\lambda_{n-1}=\lambda_{p-1}=0$. As in the proof of
Theorem \ref{final} this means that $\eta^m\neq 1$, $\eta^n\neq
1$, $\eta^p\neq 1$, $$\begin{array}{r} m(\eta^n-1)=n(\eta^m-1)\\
\mbox{ and }\; m(\eta^p-1)=p(\eta^m-1).
\end{array}$$
As before we identify ${\Bbb Q}(\eta)$ with a subfield of
${\Bbb C}$. By the previous Lemma $a=|\eta|\neq 1$. Now consider
the function $h(x)=m(a^x-1)-x(a^m-1)$. Note that
$h(m)=h(n)=h(p)=0$. We obtain a contradiction since $h'(x)$ has at
most one zero.
\end{proof}

\section{Concluding Remarks}

\subsection{\label{C1}} We apply the methods developed earlier in this paper
 to the case where $\be=0$. Let $A=A(\al,0,\gam)$.
By \cite[Lemma 4.3]{KMP} $A$ is not Noetherian. We show that there
is usually a proper homomorphic image of $A$ which is not
Noetherian. Let
               $$ w_1=du-\al ud-\gam $$
               $$ w_2=(\al -1) du+\gam. $$
Then $dw_1=w_1u=0$. Also $dw_2=\al w_2 d$ and $w_2 u=\al uw_2$. In
particular if $\al\neq 0$ then $Aw_2=w_2A$ and $A/(w_2)$ is
isomorphic to the algebra $B=k[x,y]$ generated by $x,y$ subject to
the relation $xy=\gam$. It is well known that $B$ is not
Noetherian if $\gam=0$. If $\gam\neq 0$ $B$ is not von Neuman
finite and in particular $B$ is not Noetherian. Further results on
down-up algebras $A(\al,\be,\gam)$ with $\be=0$ can be found in
\cite{Jordan98} and \cite{KK}.

\subsection{\label{C2}} We conclude with some remarks about homogenizations of down-up algebras.
Assume $\be\neq 0$. If $A=\cup A_n$ is a filtered algebra, the Rees algebra
or homogenization of the filtration is the subalgebra $\oplus A_nT^n$ of $A[T]$.

A natural way to define a filtration on $A$ is to take $A_0=K$,
$A_1$ a subspace of $A$, containing $A_0$ and a set of algebra
generators for $A$, and set $A_n=(A_1)^n$ for $n\in {\Bbb N}$.
When $A=A(\al,\be,\gam)$ is a down-up algebra an obvious choice
for $A_1$ is $A_1=span\{1,d,u\}$. We denote the Rees algebra of
the filtration obtained in this way by $H_1=H_1(\al,\be,\gam)$.
Clearly $H_1$ is generated by $D=dT$, $U=uT$ and the central
element $T$. Moreover we have
\begin{eqnarray*}
(R5)\qquad D^2U & = & d^2uT^3=(\al dud+\be ud^2+\gam d)T^3\\
                & = & \al DUD+\be UD^2+\gam DT^2.
\end{eqnarray*}

This relation is a homogenization of relation $(R1)$. Similarly $H_1$ satisfies a
homogenization of relation $(R2)$.
Thus $H_1$ is the algebra refered to in \cite[Question f)]{BenkartRoby}.
It is possible to show that $H_1$ has Hilbert series $[(1-t)^3(1-t^2)]^{-1}$,
\cite[Proposition 4.2.8]{Bauwens}.

However there is another set of generators for $A$ which resembles
the usual set of generators for $U(sl_2)$. Let $\lambda,\mu$ be
the roots of the equation $x^2-\al x -\be=0$. Thus
$\lambda+\mu=\al$, $\lambda\mu=-\be$ and set $$(R6)\qquad
h=du-\lambda ud.$$
Then
$$(R7)\qquad d h-\mu hd =\gam d$$
and
$$(R8)\qquad hu-\mu uh =\gam u.$$

By modifying the argument given in \cite[\S 3.3]{KMP}, for the case $\gam=0$, we see
that $A$ is generated by $h,u,d$ with relations $(R6), (R7), (R8)$.

Now set $A_1'=span\{1,h,u,d\}$ and let $A'_n=(A'_1)^n$ to obtain a second filtration on $A$.
The Rees algebra of this filtration is denoted $H_2=H_2(\al,\be,\gam)=\oplus A_n'T^n$.
The Hilbert series for $H_2$ is the same as that of a polynomial algebra in four variables.

\begin{thm*}
The algebras $H_1$ and $H_2$ are Auslander-regular and Cohen-Macaulay with global dimension four.
\end{thm*}
\begin{proof}
Note that $H_1/(T)\iso A(\al,\be,0)$ is Auslander regular of global dimension 3
and Cohen-Macaulay by \cite[Lemma 4.2 and Theorem 4.1]{KMP}. Using the Lemma in \cite{LS}
and writing $H_2/(T)$ as an iterated Ore extension it follows that $H_2/(T)$ also satisfies
these properties. From now on let $H$ denote either $H_1$ or $H_2$. Note that $H$ is
graded and $T$ is a homogeneous central element of positive degree in $H$ which is not
a zero divisor. Thus we can use a graded version of Nakayama's lemma and the proof of
\cite[Theorem 7.3.7]{M&R} to show that ${\gldim}(H)=4$. The result now follows from
\cite[Theorem 3.6]{L}.
\end{proof}

\begin{note*}
Parts of the theorem have been obtained independently by Bauwens
\cite[Remark 4.2.9 and Proposition 4.4.1]{Bauwens}.
The noncommutative algebraic geometry arising from the graded algebras $H_1$, $H_2$ is
studied in detail in \cite{Bauwens}. In particular the point and line modules
are obtained in the ``generic'' case.
\end{note*}

Finally we note that when $A=A(\al,\be,0)$, $A$ is a graded algebra which is
Auslander-regular of global dimension 3.
In \cite{ATV}, the regular algebras with 2 generators and 2 defining relations of degree 3
are classified in terms of a divisor $E$ in ${\Bbb P}^1\times {\Bbb P}^1$, and an automorphism
$\sigma$ of $E$. It is easily checked that
$E={\Bbb P}^1\times {\Bbb P}^1$ when $\al=0$. If $\al\neq 0$, $A=A(\al,\be,0)$ is an
algebra of type $S_1$ in \cite[4.13]{ATV}, see also \cite[Table 3.9]{AS}, that is
$E=E_1\cup E_2$ is the union of two curves of bidegree $(1,1)$ and $\sigma$ stabilizes each component.
Furthermore we have $E_1=E_2$ if and only  the equation $x^2-\al x-\be=0$ has
multiple roots. This occurs for example for $A(2,-1,0)$ which is the enveloping algebra
of the Heisenberg Lie algebra and this case is worked out in detail in
\cite[pages 36-37]{ATV}.

\end{document}